\documentclass[a4paper,12pt]{article}
\usepackage[T2A]{fontenc}
\usepackage[utf8]{inputenc} 
\usepackage{amssymb,amsfonts}
\usepackage{amsfonts,amsmath,amssymb,amscd,latexsym}
\usepackage{srcltx}

\linethickness{0.4pt}

\newcommand{\bi}{\bibitem}
\newcommand{\nb}{\newblock}

\newcommand{\be}[1]{\begin{equation}\label{#1}}
\newcommand{\ee}{\end{equation}}

\newcommand{\la}{\langle\,}
\newcommand{\ra}{\,\rangle}
\newcommand{\ve}{\varepsilon}

\newtheorem{thm}{\quad Theorem}
\newtheorem{lm}{\quad Lemma}

\title{On the Ore condition for the group ring of R.\,Thompson's group $F$}
\author{\vspace{2ex}
V. S. Guba\thanks{This work is supported by the Russian Foundation
for Basic Research, project no. 20-01-00465.}\\
Vologda State University,\\
15 Lenin Street,\\
Vologda\\
Russia\\
160600\\
E-mail: gubavs{@}vogu35.ru}
\date{}

\begin{document}

\maketitle

\begin{abstract}

Let $R=K[G]$ be a group ring of a group $G$ over a field $K$. The Ore condition says that for any $a,b\in R$ there exist $u,v\in R$ such that $au=bv$, where $u\ne0$ or $v\ne0$. It always holds whenever $G$ is amenable. Recently it was shown that for R.\,Thompson's group $F$ the converse is also true. So the famous amenability problem for $F$ is equivalent to the question on the Ore condition for the group ring of the same group.

It is easy to see that the problem on the Ore condition for $K[F]$ is equivalent to the same property for the monoid ring $K[M]$, where $M$ is the monoid of positive elements of $F$. In this paper we reduce the problem to the case when $a$, $b$ are homogeneous elements of the same degree in the monoid ring. We study the case of degree $1$ and find solutions of the Ore equation. For the case of degree $2$, we study the case of linear combinations of monomials from $S=\{x_0^2,x_0x_1,x_0x_2,x_1^2,x_1x_2\}$. This set is not doubling, that is, there are nonempty finite subsets $X\subset M\subset F$ such that $|SX| < 2|X|$. As a consequence, the Ore condition holds for linear combinations of these monomials. We give an estimate for the degree of $u$, $v$ in the above equation.
  
The case of monomials of higher degree is open as well as the case of degree $2$ for monomials on $x_0,x_1,...,x_m$, where $m\ge3$. Recall that negative answer to any of these questions will immediately imply non-amenability of $F$.

\end{abstract}

\section*{Acknowledgements}

The author is grateful to A.\,Yu.\,Ol'shanskii, Matt Brin, Mark Sapir, Rostislav Grigorchuk, and David Kielak for helpful discussions of the results of this paper.

\section*{Introduction}

Let $M$ be a monoid given by the following infinite presentation

\be{xinf}
\la x_0,x_1,x_2,\ldots\mid x_j{x_i}=x_ix_{j+1}\ (0\le i < j)\,\ra.
\ee

It is easy to see that any word in these generators can be reduced to a word of the form $x_{i_1}\ldots x_{i_k}$, where $k\ge0$ and $0\le i_1\le\cdots\le i_k$. The rewrite system $x_jx_i\to x_ix_{j+1}$ ($j > i\ge0$) turns out to be terminating and confluent so the above normal form is unique.
\vspace{1ex}

For any elements $a,b\in M$ there exists a least common right multiple of $a$, $b$. Also it is easy to check that $M$ is a cancellative monoid. A classical Ore theorem states that for a cancellative monoid $M$ with common right multiples, there exists a natural embedding of $M$ into its group of quotients (see~\cite{Ljap} for details). Any element of this group belongs to $MM^{-1}$, and the group is given be the same presentation. We denote it by $F$.
\vspace{1ex}

This group was found by Richard J. Thompson in the 60s. We refer to the survey \cite{CFP} for details. (See also \cite{BS,Bro,BG}.) It is easy to see that for any $n\ge2$, one has $x_n=x_0^{-(n-1)}x_1x_0^{n-1}$ so the group is generated by $x_0$, $x_1$. It can be given by the following presentation with two defining relations

\be{x0-1}
\la x_0,x_1\mid x_1^{x_0^2}=x_1^{x_0x_1},x_1^{x_0^3}=x_1^{x_0^2x_1}\ra,
\ee
where $a^b\leftrightharpoons b^{-1}ab$. 

Each element of $F$ can be uniquely represented by a {\em normal form\/}, that is, an expression of the form
\be{nf}
x_{i_1}x_{i_2}\cdots x_{i_k}x_{j_l}^{-1}\cdots x_{j_2}^{-1}x_{j_1}^{-1},
\ee
where $k,l\ge0$, $0\le i_1\le i_2\le\cdots\le i_k$, $0\le j_1\le j_2
\le\cdots\le j_l$ and the following is true: if (\ref{nf}) contains both $x_i$ and $x_i^{-1}$ for some $i\ge0$, then it also contains $x_{i+1}$ or $x_{i+1}^{-1}$ (in particular, $i_k\ne j_l$).
\vspace{1ex}

An equivalent definition of $F$ can be given in the following way. Let us consider all strictly increasing continuous piecewise-linear functions from the closed unit interval onto itself. Take only those of them that are differentiable except at finitely many dyadic rational numbers and such that all slopes (derivatives) are integer powers of $2$. These
functions form a group under composition. This group is isomorphic to $F$. Another useful representation of $F$ by piecewise-linear functions can be obtained if we replace $[0,1]$ by $[0,\infty)$ in the previous definition
and impose the restriction that near infinity all functions have the form $t\mapsto t+c$, where $c$ is an integer.
\vspace{1ex}

It is known that $F$ is a {\em diagram group} over the simplest semigroup presentation ${\cal P}=\langle x=x^2\rangle$. See~\cite{GbS} for the theory of these groups. In our paper, we will use diagram presentation for the elements of $F$. This is based on non-spherical diagrams over ${\cal P}$. A detailed description can be found in~\cite{Gu04}. We will mostly use diagrams that represent elements of the monoid $M$. Such objects are called {\em positive diagrams} over $x=x^2$. Let us give a brief illustration.
\vspace{1ex}

Given a normal form of an element in $F$, it is easy to draw the corresponding non-spherical diagram,
and vice versa. The following example illustrates the diagram that corresponds to the element $g=x_0^2x_1x_6x_3^{-1}x_0^{-2}$ represented by its normal form:

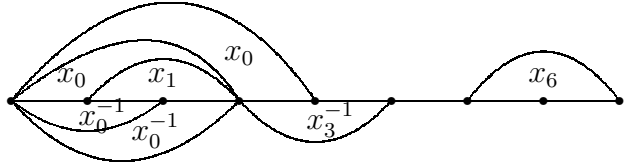
\begin{center}
	\begin{picture}(87.00,30.00)
	\put(6.00,9.00){\circle*{1.00}}
	\put(16.00,9.00){\circle*{1.00}}
	\put(26.00,9.00){\circle*{1.00}}
	\put(36.00,9.00){\circle*{1.00}}
	\put(46.00,9.00){\circle*{1.00}}
	\put(56.00,9.00){\circle*{1.00}}
	\put(66.00,9.00){\circle*{1.00}}
	\put(76.00,9.00){\circle*{1.00}}
	\put(86.00,9.00){\circle*{1.00}}
	\put(6.00,9.00){\line(1,0){80.00}}
	\bezier{132}(66.00,9.00)(76.00,22.00)(86.00,9.00)
	\bezier{120}(16.00,9.00)(26.00,20.00)(36.00,9.00)
	\bezier{176}(36.00,9.00)(26.00,25.00)(6.00,9.00)
	\bezier{264}(46.00,9.00)(28.00,35.00)(6.00,9.00)
	\bezier{120}(36.00,9.00)(46.00,-2.00)(56.00,9.00)
	\bezier{104}(6.00,9.00)(16.00,1.00)(26.00,9.00)
	\bezier{176}(6.00,9.00)(20.00,-7.00)(36.00,9.00)
	\put(36.00,15.00){\makebox(0,0)[cc]{$x_0$}}
	\put(14.00,12.00){\makebox(0,0)[cc]{$x_0$}}
	\put(26.00,12.00){\makebox(0,0)[cc]{$x_1$}}
	\put(76.00,12.00){\makebox(0,0)[cc]{$x_6$}}
	\put(25.00,5.00){\makebox(0,0)[cc]{$x_0^{-1}$}}
	\put(48.00,6.50){\makebox(0,0)[cc]{$x_3^{-1}$}}
	\put(18.00,7.00){\makebox(0,0)[cc]{$x_0^{-1}$}}
	\end{picture}
\end{center}

Here we assume that each egde of the diagram is labelled by a letter $x$. The horizontal path in the picture cuts the diagram into two parts, {\em positive} and {\em negative}. The positive part represents an element of $M$, namely, $x_0^2x_1x_6$. The top label of this positive diagram has label $x^4$, the bottom part of it has label $x^8$. So we have a positive $(x^4,x^8)$-diagram over ${\cal P}$. The set of all positive diagrams with top label $x^m$ and bottom label $x^n$, where $m\le n$, will be denoted by ${\cal S}_{m,n}$. Given a positive $(x^m,x^n)$-diagram $\Delta_1$ and a positive $(x^n,x^k)$-diagram $\Delta_2$, one can {\em concatenate} them obtaining an $(x^m,x^k)$-diagram denoted by $\Delta_1\circ\Delta_2$. This operation is very natural: we identify the bottom path of $\Delta_1$ with the top path of $\Delta_2$. 

We also say {\em product} instead of concatenation, and can multiply sets of diagrams in this way. Clearly, ${\cal S}_{m,n}{\cal S}_{n,k}={\cal S}_{m,k}$. All sets ${\cal S}_{m,n}$ are finite, and the cardinalities of them are given by numbers from Catalan triangle; see Section~\ref{oca}.
\vspace{2ex}

Recall that a group $G$ is called {\em amenable\/} whenever there exists a finitely additive normalized invariant mean on $G$, that is, a mapping $\mu\colon{\cal P}(G)\to[0,1]$ such that 

\begin{itemize}
\item 
$\mu(A\cup B)=\mu(A)+\mu(B)$ for any disjoint subsets $A,B\subseteq G$,

\item
 $\mu(G)=1$, 
 
\item
$\mu(Ag)= \mu(gA)=\mu(A)$ for any $A\subseteq G$, $g\in G$. 
\end{itemize}
One gets an equivalent definition of amenability if only one-sided invariance of the mean is assumed, say, the condition $\mu(Ag)=\mu(A)$ ($A\subseteq G$, $g\in G)$. The proof can be found in \cite{GrL}.

The class of amenable groups includes all finite groups and all abelian groups. It is invariant under taking subgroups, quotient groups, group extensions, and ascending unions of groups. The closure of the class of finite and abelian groups under these operations is the class EA of {\em elementary amenable\/} groups. A free group of rank $>1$ is not amenable. There are many useful
criteria for (non)amenability. Here we would like to mention Folner criterion from~\cite{Fol}. For our reasons, it is convenient
to formulate it as follows.
\vspace{1ex}

{\sl A group $G$ is amenable if and only if for any finite set $A\subseteq G$ and for any $\varepsilon > 0$, there exists a finite set $S$ such that $|AS| < (1+\varepsilon)|S|$.}
\vspace{1ex}

One can assume that $A$ contains the identity element. Also one can extend $A$ to the generating set of $G$ provided the group is finitely generated. In this case we see that the set $S$ from the above statement is almost invariant under (left) multiplication by elements in $A$.
\vspace{1ex}

It was proved in~\cite{BS} that the group $F$ has no free subgroups of rank $>1$. It is also known that $F$ is not elementary amenable (see~\cite{CGH}). However, the famous problem about amenability of $F$ is still open. The question whether $F$ is amenable was asked by Ross Geoghegan in 1979; see~\cite{Ger87}. There were many attempts of various authors to solve this problem in both directions. We will not review a detailed history of the problem; this information can be found in many references. However, to emphasize the difficulty of the question, we mention the paper~\cite{Moore13}, where it was shown that if $F$ is amenable, then Folner sets for it have a very fast growth. Besides, we would like to refer to the paper~\cite{BB05} where the authors obtained an estimate of the isoperimetric constant of the group $F$ in its standard set of generators $\{x_0,x_1\}$. This estimate has not been improved so far.
\vspace{1ex}

Now we are going to describe a new approach to the problem in terms of equations in the group ring of $F$. 

\section{Ore condition and amenability}
\label{oca}

Tamari~\cite{Ta54} shows that if a group $G$ is amenable, then the group ring $R=K[G]$ satsfies Ore condition for
any field $K$. This means that for any $a,b\in R$ there exist $u,v\in R$ such that $au=bv$, where $u\ne0$ or $v\ne0$. 

Let us generalize this statement. Suppose that instead of one linear equation $au=bv$ with coefficients in $R$ we have a system of them, where the number of variables exceeds the number of equations:
$$
\left\{
\begin{array}{ccccccc}
a_{11}u_1&+  &\cdots  &+  &a_{1n}u_n  & = & 0 \\
\cdots&  &  \cdots & &\cdots  &  &  \\
a_{m1}u_1&+  &\cdots  &+  &a_{mn}u_n  &=  &0 
\end{array}
\right.
$$
where $n > m$, $a_{ij}\in R$ for all $1\le i\le m$, $1\le j\le n$. We are interested in solutions $(u_1,...,u_n)\in R^n$. 

We claim that for amenable group $G$, this system always has a nonzero solution.
\vspace{1ex}

Indeed, let $A\subseteq G$ be the union of supports of all the coefficients of the form $a_{ij}$. By F\o{}lner criterion, for any $\varepsilon > 0$ there exists a finite set $X$ such that $|AX| < (1+\varepsilon)|X|$. Let $u_j$ ($1\le j\le n$) be linear combinations of elements in $X$ with indefinite coefficients from the field. This gives $n|X|$ variables taking their values in $K$. For any of $m$ equations in the system, we collect all terms on any element in $AX$ and impose the condition that the sum of them is zero. This leads to an ordinary system of $m|AX|$ linear equations. Such a system has a nonzero solution whenever the number of variables exceeds the number of equations. This holds if $m(1+\varepsilon) < n$, so it suffices to claim $\varepsilon < \frac{n-m}m$.
\vspace{1ex}

In a recent paper~\cite{Ba19}, Bartholdi shows that the converse to the above statement is true. This gives a new criterion for amenabilty of groups. Although Thorem 1.1 in~\cite{Ba19} concerns the so-called GOE and MEP properties of automata (Gardens of Eden and Mutually Erasable Patterns), the proof of it allows one to extract the following statement.

\begin{thm}
\label{barth19}
{\rm(Bartholdi)}\ For any group $G$, the following two properties are equaivalent.

(i) $G$ is amenable

(ii) For any field $K$ and for any system of $m$ linear equations over $R=K[G]$ in $n > m$ variables, there exists a nonzero solution.
\end{thm}

In the Appendix to the same paper, Kielak shows that if the group ring $K[G]$ has no zero divisors, both properties are equivalent to the Ore condition. In particular, this holds for R.\,Thompson's group $F$. It is orderable, so there are no zero divisors in a group ring over a field. So we quote the following

\begin{thm}
\label{kielak}
{\rm(Kielak)}\ The group $F$ is amenable if and only if the group ring $K[F]$ over any field satisfies Ore condition.
\end{thm}

The following elementary property of the group $F$ is well known.

\begin{lm}
\label{gig}
For any $g_1,...,g_n\in F$ there exist $g\in F$ such that $g_1g,\ldots,g_ng\in M$.
\end{lm}

Using this fact, one can easily show that the Ore condition for $K[F]$ is equivalent to the Ore condition for $K[M]$.

\begin{lm}
\label{kfkm}
For any field $K$, the group ring $K[F]$ satisfies Ore condition if and only if the monoid ring $K[M]$ satisfies Ore condition.
\end{lm}
\vspace{1ex}

\begin{prf}
Suppose that any equation of the form $au=bv$, where $a,b\in K[F]$, has a nonzero solution in the same ring. Let $a,b\in K[M]$. Then there exist $u,v\in K[F]$, where $u\ne0$ or $v\ne0$, such that $au=bv$. Let $X$ be the union of supports of $u$ and $v$. This is a finite subset in $F$. According to Lemma~\ref{gig}, there exists $g\in F$ such that $Xg\subset M$. Therefore, $aug=bvg$, where $(ug,vg)$ is a non-zero solution of the equation in $K[M]$.

Converesly, suppose that the Ore condition holds in $K[M]$. Take any $a,b\in K[F]$. Taking $Y$ as a union of their supports, find $g\in G$ such that $Yg\subset M$. Then $ag$, $bg$ belong to $K[M]$ and so there exists a nonzero solution of the equation $ag\cdot u=bg\cdot v$ in $K[M]$. Here $(gu,gv)$ brings a nonzero solution of the equation in $K[F]$ with coefficients $a$, $b$.
\end{prf}
\vspace{1ex}

The role of this easy lemma is that the monoid ring $K[M]$ has a simpler structure than the group ring of $F$. This is in fact a graded ring of skew polynomials where the non-commutative variables are multiplied using the rule $x_jx_i=x_ix_{j+1}$ ($j > i$).

Now we are going to reduce the general problem of solving equations to the one for homogeneous polynomials $a$, $b$.

\begin{lm}
\label{homogen}
Suppose that any equation of the form $au=bv$ has a nonzero solution in $K[M]$ provided $a$, $b$ are homogeneous polynomials of the same degree. Then $K[M]$ satisfies Ore condition.
\end{lm}

\begin{prf}
Notice that if both $a$, $b$ are homogeneous but not of the same degree then one can multiply one of them on the right by a monomial such that the degrees will coincide. Say, if $m=\deg a < \deg b=n$ then we solve the equation $ax_0^{n-m}u=bv$ for monomials of the same degree, and this gets us the solution for the pair $a,b$.

Also it is easy to see that if $au=bv$, then one can replace $u$, $v$ by their homogeneous components of minimal degree. This will give another solution to the same equation where $u$, $v$ are also homogeneous.

Any nonzero element $a\in K[M]$ we can written as a sum of homogeneous components: $a=a_1+\cdots+a_k$, where $k\ge1$, $\deg a_1 < \cdots < \deg a_k$. By the {\em width\/} of $a$ we mean the difference of highest and lowest degrees: $\alpha=\deg a_k-\deg a_1$. Similarly, let $b=b_1+\cdots+b_l$ be a sum of its homogeneous components, where $l\ge1$, $\deg b_1 < \cdots < \deg b_l$, and the width of $b$ is $\beta=\deg b_l-\deg b_1$.

We proceed by induction on $\alpha+\beta$. The cases $a=0$ and $b=0$ are obvious so we will assume that both polynomials $a$, $b$ are nonzero. The situation $\alpha+\beta=0$ means that $a$, $b$ are homogeneous. In this case the equation has a nonzero solution.

Let the sum of widths be greater than zero, and assume that $\alpha\le\beta$ without loss of generality. Here $\beta > 0$ so $l\ge2$. Let us solve the equation $a_1u_1=b_1v_1$ for homogeneous $a_1$, $b_1$, where $u_1$, $v_1$ are homogeneous and nonzero. 

Let $b'=bv_1-au_1=b_2v_1+\cdots+b_lv_1-a_2u_1-\cdots-a_ku_1$. The terms $a_1u_1$ and $b_1v_1$ cancel in this difference. The case $b'=0$ already gives us a nonzero solution of the equation for $a$, $b$. So assume that $b'\ne0$.

Let us compare the homogeneous components in $b'$ of highest and lowest degree. The highest one does not exceed $\max(\deg(b_lv_1),\deg(a_ku_1))$. We see that $\deg(b_lv_1)=\deg b_l+\deg v_1=\beta+\deg b_1v_1$ and $\deg(a_ku_1)=\deg a_k+\deg u_1=\alpha+\deg(a_1u_1)$. Since $a_1u_1$ and $b_1v_1$ are equal, we have $\deg(b_lv_1)-\deg(a_ku_1)=\beta-\alpha\ge0$. Therefore, the highest degree of homogeneous components in $b'$ does not exceed $\deg b_l+\deg v_1$.

On the other hand, all monomials in $b'$ have degree strictly greater that $\deg b_1+\deg v_1=\deg a_1+\deg u_1$ since degrees of the $a_i$s and the $b_j$s strictly increase. So the difference between the highest and the lowest degree of monomials in $b'$ is strictly less than $(\deg b_l+\deg v_1)-(\deg b_1+\deg v_1)=\beta$. Hence the sum of widths of $a$ and $b'$ is strictly less than $\alpha+\beta$, and we can use the inductive assumption.

Solving the equation $au=b'v$, we find its nonzero solution (both $u$, $v$ are nonzero). Now $au=(bv_1-au_1)v$ so we have $a(u+u_1v)=bv_1v$. This solution is nonzero because $v_1,v\ne0$, and the group ring $K[F]\supset K[M]$ has no zero divisors. This completes the proof.
\end{prf}
\vspace{1ex}

Now we have a bunch of equations in $K[M]$ indexed by two parameters. The fist one is $d$, the degree of homogeneous polynomials $a$ and $b$. The second one is $m$, where $m$ is the highest subscript in variables we involve. The general strategy can be the following: we try to solve as much equations in $K[M]$ as we can, using this classification. For a pair of numbers $d\ge1$, $m\ge1$, we can take $a$, $b$ as linear combinations of monomials of degree $d$ in variables $x_0$, $x_1$, ... , $x_m$ with indefinite coefficients. We can think about these coefficients as elements of the field of rational functions over $K$ with a number of variables.
\vspace{1ex}

More precisely, we can state the general problem as follows. Any finite system of monomials of degree $d\ge1$ is contained in a set of the form ${\cal S}_{m+1,m+d+1}$ for some $m\ge1$. This set consists of all elements in the monoid $M$ with normal forms $x_{i_1}\ldots x_{i_d}$, where $i_1\le\cdots\le i_d$ and $i_1\le m$, $i_2\le m+1$, ... , $i_d\le m+d-1$. Let $K[S]$ denote the set of all linear combinations of elements of $S\subset M$ with coefficients in $K$.
\vspace{1ex}

{\bf Problem} ${\cal P}_{d,m}$: {\sl Given two elements $a,b\in K[{\cal S}_{m+1,m+d+1}]$, find a nonzero solution of the equation $au=bv$, where $u,v\in K[M]$, or prove that it does not exist.}
\vspace{1ex}

According to Theorem~\ref{kielak} by Kielak, and Lemma~\ref{homogen} on homogeneous equations, we have the following alternative. If the Problem ${\cal P}_{d,m}$ has positive solution for any $d,m\ge1$ (that is, we can find nonzero solutions), then the group $F$ is amenable. If this Problem has negative solution for at least one case, then $F$ is not amenable.
\vspace{1ex}

\section{Cardinality arguments}
\label{card}

The first case we are going to start with, is the case of polynomials of degree $d=1$ for arbitrary $m$. Let us show that equations of the form
\be{ab}
(\alpha_0x_0+\alpha_1x_1+\cdots+\alpha_mx_m)u=(\beta_0x_0+\beta_1x_1+\cdots+\beta_mx_m)v
\ee
have nonzero solutions in $K[M]$ for arbitrary coefficients $\alpha_i,\beta_i\in K$ ($0\le i\le m$). Our first approach will be based on cardinality reasons; later we are going to improve it.

All elements of the set $X_m=\{x_0,x_1,...,x_m\}$ are viewed as positive diagrams over $x=x^2$ with top label $x^{m+1}$ and bottom label $x^{m+2}$. Each of these diagrams consists of one cell; the diagram for $x_i$ ($0\le i\le m$) looks in the following way:

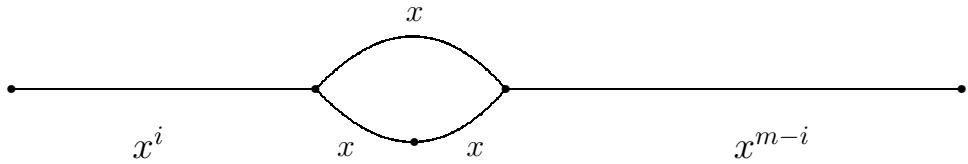
\begin{center}
	\begin{picture}(131.00,25.00)
	\put(2.00,11.00){\circle*{1.00}}
	\put(2.00,11.00){\line(1,0){40.00}}
	\put(42.00,11.00){\circle*{1.00}}
	\put(67.00,11.00){\circle*{1.00}}
	\bezier{152}(42.00,11.00)(55.00,25.00)(67.00,11.00)
	\bezier{152}(42.00,11.00)(55.00,-3.00)(67.00,11.00)
	\put(67.00,11.00){\line(1,0){60.00}}
	\put(55.00,4.00){\circle*{1.00}}
	\put(127.00,11.00){\circle*{1.00}}
	\put(55.00,21.00){\makebox(0,0)[cc]{$x$}}
	\put(46.00,3.00){\makebox(0,0)[cc]{$x$}}
	\put(63.00,3.00){\makebox(0,0)[cc]{$x$}}
	\put(20.00,4.00){\makebox(0,0)[cc]{\large$x^i$}}
	\put(102.00,4.00){\makebox(0,0)[cc]{\large$x^{m-i}$}}
	\end{picture}
\end{center}

We will need some more notation. An edge labelled by $x^i$ will be denoted $\varepsilon(x^i)$. Also if we have two diagrams $\Delta_1$ and $\Delta_2$, then $\Delta_1+\Delta_2$ denotes the diagram where the rightmost point of $\Delta_1$ is identified with the leftmost point of $\Delta_2$. The diagram in the picture can thus be denoted as $\varepsilon(x^i)+(x=x^2)+\varepsilon(x^{m-i})$, where $(x=x^2)$ is a positive cell.

All these diagrams belong to the set ${\cal S}_{m+1,m+2}$. Let $n\gg1$ be a sufficiently large integer. We are going to multiple the set $X_m$ that we identify with ${\cal S}_{m+1,m+2}$, by the set $Y={\cal S}_{m+2,n}$. Our aim is to show that $|X_mY| < 2|Y|$. This will imply that equation~(\ref{ab}) has nonzero solutions.

Notice that the set $X_m$ has cardinality $m+1$ which can be arbitrarily large. The list of elements of $X_mY$ taken with repetitions will have size $(m+1)|Y|$. Nevertheless, it has less than $2|Y|$ distinct elements.

It is known from standard combinatorics that the cardinality of ${\cal S}_{k,n}$ is given by an element of the Catalan triangle:
\be{cat}
b_{nk}=\frac{k(2n-k-1)!}{n!(n-k)!}.
\ee
So in our case we have 
$$
\frac{|X_mY|}{|Y|}=\frac{|{\cal S}_{m+1,n}|}{|{\cal S}_{m+2,n}|}=\frac{b_{n,m+1}}{b_{n,m+2}}=\frac{m+1}{m+2}\cdot\frac{2n-m-2}{n-m-1}.
$$
The desired inequality holds whenever $n > \frac{(m+1)(m+2)}2$. Here the cardinality approach gives us a quadratic estimate to the degree of $u$ and $v$ in a potential solution. In fact, we can improve this estimate finding solutions of degree just $m$ in variables $x_0$, $x_1$, ... , $x_{2m}$, which we are going to present in the next Section finding solutions of potentially minimal degree for~(\ref{ab}). Now let us finish our proof. 

We need to repeat the same argument that was used in the proof of Tamari's theorem in the beginning. Namely, we take $u$, $v$ as linear combinations of the elements in $Y$ with free coefficients (variables). The total number of variables is $2|Y|$. Making a substitution of these expressions into~(\ref{ab}), we get a linear combination of $|X_mY|$ elements of $M$. Claiming that all of them have equal coefficients in the left-hand side and the right-hand side, we get $|X_mY| < 2|Y|$ linear equations over $K$. Since the system with this property has a nonzero solution in the field, we get the desired solution for our equation~(\ref{ab}).

So let us state a weak form of our result on the equations for the case of polynomials of degree $1$ as follows.

\begin{thm}
\label{xmy}
a) For any $m\ge1$, the set of elements $X_m=\{x_0,x_1,...,x_m\}$ is not doubling, that is, there exists a finite subset $Y\subset M$ such that $|X_mY| < 2|Y|$. 

b) If $a,b\in K[M]$ are linear combinations of monomials $x_0$, $x_1$, ... , $x_m$ of degree $1$, then the equation $au=bv$ in $k[M]$ has a nonzero solution, where $\deg u=\deg v\le\frac{m(m+1)}2$.
\end{thm}

Here the degrees of $u$ and $v$ do not exceed $n-(m+2)$, so the above estimate holds if we take $n=\frac{(m+1)(m+2)}2+1$.
\vspace{1ex}

In other words, Problem ${\cal P}_{1,m}$ has positive solution for any $m\ge1$.
\vspace{1ex}

Now let $d=2$, that is, we deal with linear combinations of monomials of degree two. The first interesting case is $m=1$. This means that $a$, $b$ belong to $K[S]$, where $S={\cal S}_{m+1,m+d+1}={\cal S}_{2,4}=\{x_0^2,x_0x_1,x_0x_2,x_1^2,x_2^2\}$. The above is just the set of all positive $(x^2,x^4)$-diagrams over $x=x^2$ (each of them has two cells). We are going to show that for this case the Problem ${\cal P}_{d,m}={\cal P}_{2,1}$ has positive solution.
\vspace{1ex}

We will use the construction from~\cite{Don11}. Notice that the set $S_{03}$ in the notation of that paper is what we call ${\cal S}_{1,4}$, if to replace forests by diagrams. As a set of elements of $M$, this is exactly $x_0\{x_0^2,x_0x_1,x_0x_2,x_1^2,x_2^2\}$. The first $x_0$ is ignored in our case. What is called {\em ruinous} in Donnelly's paper, is called {\em doubling} in our terminilogy.

So let $S=\{x_0^2,x_0x_1,x_0x_2,x_1^2,x_2^2\}={\cal S}_{2,4}$. If one takes $Y={\cal S}_{4,n}$, then the condition $|SY| < 2|Y|$ never holds. The idea of the proof of~\cite[Theorem 1]{Don11} is as follows (using our notation). Take $n\gg1$ and exclude from ${\cal S}_{4,n}$ all diagrams of the form $\ve(x)+\Delta+\ve(x)+\ve(x)$ and $\ve(x)+\ve(x)+\Delta+\ve(x)$. Here $\Delta$  is any $(x,x^{n-3})$-diagram. Let $Y$ be the corresponding set difference; its cardinality is $|Y|=b_{n4}-2c_{n-4}$, where $c_k$ is the $k$th Catalan number, and $b_{n4}=\frac{4(2n-5)!}{n!(n-4)!}$ is the number from Catalan triangle.
\vspace{1ex}

The product $SY$ is contained in ${\cal S}_{2,n}$. It is not hard to describe the elements that do not appear in the difference of these sets. To do that, we introduce a few technical concepts.

Let $\Xi$ be a positive semigroup diagram over $x=x^2$. We call it {\em simple} whenever it has exactly one top cell (that is, the cell whose top boundary is contained in the top path of $\Xi$). In the language of rooted binary forests, this is equivalent to the property that exactly one tree of the forest is nontrivial. Suppose that we removed the top cell of $\Xi$. If the result is again a simple diagram, then we call $\Xi$ a {\em 2-simple} diagram.

Now look at those 2-simple diagrams from ${\cal S}_{2,n}$ such that if we remove their top cells two times, we get a diagram of the form  $\ve(x)+\Delta+\ve(x)+\ve(x)$ or $\ve(x)+\ve(x)+\Delta+\ve(x)$. Now it is easy to check that there are 3 ways to add top cells twice for any of these diagrams in order to get a 2-simple preimage. We hope the reader can check this geometric fact drawing a few pictures. Hence the set $SY$ differs from ${\cal S}_{24}$ by at least $6c_{n-4}$ elements. Therefore we get the following estimate:
$$
\frac{|SY|}{|Y|}\le\frac{b_{n2}-6c_{n-4}}{b_{n4}-2c_{n-4}}=\frac13\cdot\frac{29n^3-231n^2+562n-420}{5n^3-47n^2+142n-140}.
$$
Clearly, this quotient is less than 2 for $n\gg1$. In fact, the inequality $|SY| < 2|Y|$ holds for all $n\ge45$. We proved

\begin{thm}
\label{s24}
The equation of the form
$$
(\alpha_{00}x_0^2+\alpha_{01}x_0x_1+\alpha_{02}x_0x_2+\alpha_{11}x_1^2+\alpha_{12}x_1x_2)u=
(\beta_{00}x_0^2+\beta_{01}x_0x_1+\beta_{02}x_0x_2+\beta_{12}x_1^2+\beta_{12}x_1x_2)v
$$
in the monoid ring $K[M]$ has a nonzero solution with the property $\deg u,v\le41$.
\end{thm}

Here $\alpha_{ij}$, $\beta_{ij}$ were arbitrary coefficients from $K$. The equation has a solution in $K[Y]$, where $Y\subset{\cal S}_{4n}$ consist of monomials of degree $n-4\le41$. Notice that we do not know what is the minimum degree of $u$, $v$ for solutions of this equation. A rough computer search can show that $n > 10$, but we do not even know whether the minimum value of $n$ is close to $10$ or $40$. 

We think it is interesting to find the value of $n$ by the following reasons. If we go towards amenability of $F$, we need to be able to solve more and more complicated equations in $K[M]$. We can do that for $d=1$ and arbitrary $m$; we also know the answer for the case $d=2$, $m=1$. The cases that come after that are already unknown. This is $d=2$, $m=2$, where $a$, $b$ are linear combinations of 9 monomials of degree 2:
$$
x_0^2,\, x_0x_1,\, x_0x_2,\, x_0x_3,\, x_1^2,\, x_1x_2,\, x_1x_3,\, x_2^2,\, x_2x_3.
$$
This is one possible candidate to obtain the negative answer. If, nevertheless, this Problem ${\cal P}_{2,2}$ has positive answer (that is, there exists nonzero solutions), then one can try Problem ${\cal P}_{3,1}$, where $a$, $b$ are linear combinations of 14 monomials of degree 3:
$$
x_0^3, x_0^2x_1, x_0^2x_2, x_0^2x_3, x_0x_1^2, x_0x_1x_2, x_0x_1x_3, x_0x_2^2, x_0x_2x_3, x_1^3, x_1^2x_2, x_1^2x_3,x_1x_2^2,x_1x_2x_3.
$$

Cardinality reasons we used above are not so strong tools. So the strategy has to be as follows: starting from equations of a simple form, we try not only to prove they have nonzero solutions (which is not hard), but also try to describe somehow the set of all their solutions.

Say, if we have equation of the form $au=bv$, where $a=\alpha_0x_0+\alpha_1x_1$, $b=\beta_0x_0+\beta_1x_1$, then the description of all its solutions is easy. Namely, $u=(\beta_0x_0+\beta_1x_2)w$, $v=(\alpha_0x_0+\alpha_1x_2)w$ for any $w\in K[M]$. This means that the intersection of two principal right ideals $aR\cap bR$ is a principal right ideal ($R=K[M]$).

We also know how to describe all solutions of the equation $(\alpha_0x_0+\alpha_1x_1+\alpha_2x_2)u=(\beta_0x_0+\beta_1x_1+\beta_2x_2)v$. This will be done in the next Section. For this equation, the description has more complicated form; the intersection $aR\cap bR$ is no longer a principal right ideal. 
\vspace{1ex}

If instead of one equation we have a system of equations of the form
$$
(\alpha_1x_0+\beta_1x_1)u_1=(\alpha_2x_0+\beta_2x_1)u_2=\cdots=(\alpha_kx_0+\beta_kx_1)u_k
$$
for any $k\ge2$, then it also has a nonzero solution. This obviously follows from cardinality reasons: one can take a set of the form $Y=\{1,x_0^{-1}x_1,(x_0^{-1}x_1)^2,\ldots,(x_0^{-1}x_1)^n\}$ for $n\gg1$ and notice that the cardinality of $Y$ almost coincides with the cardinality of $\{x_0,x_1\}Y$. This construction is quite trivial, so we will offer a more explicit form of the solution. The product $(\alpha_1x_0+\beta_1x_1)(\alpha_1x_0+\beta_1x_2)\cdots(\alpha_kx_0+\beta_kx_{k+1})$ is left divisible by  $\alpha_ix_0+\beta_ix_1$ for any $1\le i\le k$, which can be checked directly.
\vspace{1ex}

A much more interesting example of a system of equations looks as follows. Let us state it as a separate problem.
\vspace{1ex}

{\bf Problem} ${\cal Q}_k$: {\sl Given $k+1$ linear combinations of elements $x_0$, $x_1$, $x_2$, consider a system of $k$ equations with $k+1$ unknowns:
$$
(\alpha_0x_0+\beta_0x_1+\gamma_0x_2)u_0=(\alpha_1x_0+\beta_1x_1+\gamma_1x_2)u_1=\cdots=(\alpha_kx_0+\beta_kx_1+\gamma_kx_2)u_k.
$$
Find a nonzero solution of this system, where $u_0,u_1,\ldots,u_k\in K[M]$, or prove that it does not exist.}
\vspace{1ex}

Notice that ${\cal Q}_1$ has been already considered. To solve ${\cal Q}_2$ in positive, it suffices to find a finite set $Y$ with the property $|AY| < \frac32|Y|$, where $A=\{x_0,x_1,x_2\}$. This can be done by cardinality reasons simliar to the above proof of Theorem~\ref{s24}. The estimate there will be also $n\ge45$. We do not give details here since we are able to prove a much stronger fact. Namely, using the result of~\cite{Gu19}, we can construct a finite set $Y$ with the property $|AY| < \frac43|Y|$. The size of $Y$ is really huge, it does not have transparent description. This immediately implies that Problem ${\cal Q}_3$ has a positive solution. We announce this result here; details will appear in forthcoming papers.
\vspace{1ex}

Donnelly shows in~\cite{Don14} that $F$ is non-amenable if and only if there exists $\varepsilon > 0$ such that for any finite set $Y\subset F$, one has $|AY|\ge(1+\varepsilon)|Y|$, where $A=\{x_0,x_1,x_2\}$ (see also~\cite{Don07}). For the set $Y$ here, one can assume without loss of generality that $Y$ is contained in ${\cal S}_{4,n}$ for some $n$. This gives some evidence that the amenability problem for $F$ has very close relationship with the family of Problems ${\cal Q}_k$. The case $k=4$ looks as a possible candidiate to a negative solution (that is, all solutions are zero). If true, this will imply that the constant $\varepsilon=\frac14$ fits into the above condition.
\vspace{1ex}

Descriptions of solutions for equations or their systems seem also be important because if we try to prove that some equation or system has only zero solutions, then we are based on the description of simpler equations for which we know all solutions. This situation is vaguely similar to the classical case when the descrition of all Pythagorian triples implies Fermat's proof that $X^4+Y^4=Z^2$ has no solutions in positive integers. Here we can expect a similar effect.

\section{Solutions of Small Degree}
\label{msol}

In this Section we go back to the equation~(\ref{ab}) finding its solutions of possibly minimal degree. Namely, instead of quadratic power of $\deg u$, $\deg v$ with respect to $m$, we will find solutions with $\deg u=\deg v=m$. It looks truthful that the value $m$ here is minimal. We do not offer the proof of this fact although we give a complete description to the set of solutions of our equations for the cases $m=1$ and $m=2$.

For our needs it will be convenient to assume that the coefficients in the equation are independent variables, that is, they belong to the field of rational functions over some field $L$. By $K$ we denote the field of rational functions with coefficients in $L$ over a number of variables of the form $\alpha_i$, $\beta_i$ ($i\ge0$).

So let $R=K[M]$ and let
$$
a=\alpha_0x_0+\alpha_1x_1+\cdots+\alpha_mx_m,\ \ b=\beta_0x_0+\beta_1x_1+\cdots+\beta_mx_m.
$$
We are interested in finding nonzero solution of the equation $au=bv$, where $u$, $v$ are homogeneous polynomials in $R$ of possibly minimal degree.

\begin{thm}
\label{aubv}
The equation
$$
(\alpha_0x_0+\alpha_1x_1+\cdots+\alpha_mx_m)u=(\beta_0x_0+\beta_1x_1+\cdots+\beta_mx_m)v
$$
has a nonzero solution, where $u$, $v$ are homogeneous polynomials of degree $m$ in variables $x_0,x_1,\ldots,x_{2m}$.
\end{thm}

\begin{prf}
We keep the notaion $a$, $b$ for the coefficients of the above equation. We intoduce a sequence of elements $\sigma_i=a\cdot\beta_i+b\cdot(-\alpha_i)\in aR+bR$ for all $0\le i\le m$, where the coefficient on $x_i$ in $\sigma_i$ is zero. So each $\sigma_i$ has the form $$\gamma_{i0}x_0+\cdots+\gamma_{i,i-1}x_{i-1}+\gamma_{i,i+1}x_{i+1}+\cdots+\gamma_{im}x_m$$ for any $0\le i\le m$, where all $\gamma_{ij}=\alpha_j\beta_i-\alpha_i\beta_j$ are nonzero for $j\ne i$.

We introduce polynomials of the form $f_k(x_k,\ldots,x_{k+m-1})$ of degree $k$ for each $1\le k\le m$. Each of these polynomials will be equal to $au_{k-1}+bv_{k-1}$ for some homogeneous polynomials of degree $k-1$. This will be done by induction. At each step we will check that monomial $x_0^{k-1}$ occurs in both $u_{k-1}$, $v_{k-1}$ with nonzero coefficient.

For the case $k=1$, we take $f_1(x_1,\ldots,x_m)=\sigma_0$. By definition, this is a linear combination of $x_1$, ... , $x_m$ with nonzero coefficients. It also equals $au_0+bv_0$, where $u_0=\beta_0$, $v_0=-\alpha_0$. Notice that coefficients on $x_0^0$ here are nonzero.

Let $k < m$. Multiplying $f_k(x_k,\ldots,x_{k+m-1})$ by $\sum\limits_{i=0}^{k-1}\gamma_{ki}x_i$ on the right and using defining relations of the monoid $M$, we get the element
$$
f_k(x_k,\ldots,x_{k+m-1})\cdot\sum\limits_{i=0}^{k-1}\gamma_{ki}x_i=\sum\limits_{i=0}^{k-1}\gamma_{ki}x_i\cdot f_k(x_{k+1},\ldots,x_{k+m}).
$$
We used the fact that for any $i < k$ and for any product $w$ of $x_k$, $x_{k+1}$, ... , it holds the equality $wx_i=x_iw'$, where $w'$ is obtained from $w$ by increasing all subscripts by $1$. This follows directly from the defining relations of $M$. 

Notice that $\sigma_k=\sum\limits_{i=0}^{k-1}\gamma_{ki}x_i+\sum\limits_{i=k+1}^{m}\gamma_{ki}x_i$. Therefore, we can rewrite the above equality as follows:
$$
(au_{k-1}+bv_{k-1})\cdot\sum\limits_{i=0}^{k-1}\gamma_{ki}x_i=(\sigma_k-\sum\limits_{i=k+1}^{m}\gamma_{ki}x_i)\cdot f_k(x_{k+1},\ldots,x_{k+m}).
$$
Now let us take into account that $\sigma_k=a\cdot\beta_k+b\cdot(-\alpha_k)$. Also let us define the polynomial
$$
f_{k+1}(x_{k+1},\ldots,x_{k+m})=\sum\limits_{i=k+1}^{m}\gamma_{ki}x_i\cdot f_k(x_{k+1},\ldots,x_{k+m}).
$$
This is a homogeneous polynomial of degree $k+1$ in variables $x_{k+1}$, ... , $x_{k+m}$. Sending this polynomial to the left-hand side of the above equality, we see that it is equal to
$$
-(au_{k-1}+bv_{k-1})\cdot\sum\limits_{i=0}^{k-1}\gamma_{ki}x_i+(a\cdot\beta_k+b\cdot(-\alpha_k))\cdot f_k(x_{k+1},\ldots,x_{k+m}).
$$
This means that $f_{k+1}(x_{k+1},\ldots,x_{k+m})=au_k+bv_k$, where
$$
u_k=-u_{k-1}\sum\limits_{i=0}^{k-1}\gamma_{ki}x_i+\beta_kf_k(x_{k+1},\ldots,x_{k+m}),
$$
$$
v_k=-v_{k-1}\sum\limits_{i=0}^{k-1}\gamma_{ki}x_i-\alpha_kf_k(x_{k+1},\ldots,x_{k+m}).
$$
These are homogeneous polynomials of degree $k$. Notice that the term $x_0^k$ occurs in both expressions with nonzero coefficient. Indeed, $x_0^{k-1}$ had nonzero coefficients in $u_{k-1}$, $v_{k-1}$; here it is multiplied by $-\gamma_{k0}x_0$, and the rest of the sum in both expressions does not contain $x_0$ at all.
\vspace{1ex}

Continuing in this way, at the last step we get $f_m(x_m,...,x_{2m-1})$ as a homogeneous polynomial of degree $m$. It has the form $au_{m-1}+bv_{m-1}$. Now we multiply it on the right by $\sigma_m=\gamma_{m0}x_0+\cdots+\gamma_{m,m-1}x_{m-1}$. Aplying defining relations of $M$, we get to
$$
f_m(x_m,...,x_{2m-1})(\gamma_{m0}x_0+\cdots+\gamma_{m,m-1}x_{m-1})=(\gamma_{m0}x_0+\cdots+\gamma_{m,m-1}x_{m-1})f_m(x_{m+1},...,x_{2m}).
$$
Now $f_m\sigma_m$ becomes $au_{m-1}\sigma_m+bv_{m-1}\sigma_m$. On the other hand, it equals $\sigma_mf_m(x_{m+1},...,x_{2m})$, where $\sigma_m\in aR+bR$ is expressed according to its defintion. So the same element also equals $(a\beta_m-b\alpha_m)f_m(x_{m+1},...,x_{2m})$. 

Comparing both things, we have an equality of the form $au=bv$, where 
$$
u=-u_{m-1}\sigma_m+\beta_mf_m(x_{m+1},...,x_{2m}),\ \ \ v=v_{m-1}\sigma_m+\alpha_mf_m(x_{m+1},...,x_{2m}).
$$
By the same reasons as above, both $u$ and $v$ have monomial $x_0^m$ with nonzero coefficient (the monomials $x_0^{m-1}$ with nonzero coefficients from $u_{m-1}$, $v_{m-1}$ are multiplied by $\pm\gamma_{m0}x_0$). 

It is also clear that $u$, $v$ are homogeneous polynomials of degree $m$ in variables $x_0$, $x_1$, ... , $x_{2m}$, as it was stated. This completes the proof.
\end{prf}
\vspace{1ex}

Now let us observe the case $m=1$ in order to describe the set of all solutions. Here $a=\alpha_0x_0+\alpha_1x_1$, $b=\beta_0x_0+\beta_1x_1$ will be linear combinations of $x_0$, $x_1$ with arbitrary coefficients from any field. We exclude trivial case when the elements $a$, $b$ are proportional. The equation $au=bv$ now can be rewritten as $x_0(\alpha_0u-\beta_0v)=x_1(\beta_1v-\alpha_1u)$. By our assumption, both expressions here are nonzero. All monomials in the right-hand side are left divisible by $x_0$. Thus we can write $\beta_1v-\alpha_1u=x_0w$ for some polynomial $w\in K[M]$. Using the fact that $x_1x_0=x_0x_2$ and cancelling by $x_0$ on the left, we get $\alpha_0u-\beta_0v=x_2w$. Solving the above system, one can write $u=(\beta_0x_0+\beta_1x_2)w$, $v=(\alpha_0x_0+\alpha_1x_2)w$ without loss of generality, where $w$ was replaced by $(\alpha_0\beta_1-\alpha_1\beta_0)w$ up to a nonzero coefficient.

Now let us go to a more interesting case $m=2$ in order to give a complete description of the set of solutions for $au=bv$, where $a=\alpha_0x_0+\alpha_1x_1+\alpha_2x_2$, $b=\beta_0x_0+\beta_1x_1+\beta_2x_2$. To avoid unnecessary work with many coefficients, let us apply a linear transformation, reducing the equation to the following one: 
\be{012}
(x_0+\alpha x_2)u=(x_1+\beta x_2)v,
\ee
where $\alpha,\beta$ are some coefficients. We will assume both of them are nonzero: otherwise the description turns out to be trivial.

First of all, we take a solution of this equation extracted from the proof of Theorem~\ref{aubv}. One can check directly that the following polynomials satisfy~(\ref{012}):
$$
u_0=\beta x_0x_3+\beta^2x_0x_4-\alpha x_1x_3-\alpha\beta x_1x_4-\alpha\beta x_3^2-\alpha\beta^2 x_3x_4,
$$
$$
v_0=\beta x_0^2-\alpha x_0x_1-\alpha^2x_3^2-\alpha^2\beta x_3x_4.
$$
We say that $(u_0,v_0)$ is a {\em basic} solution. Now we are going to show how to extract all solutions from it.
\vspace{1ex}

By $M_1$ we denote the submonoid of $M$ generated by $x_1$, $x_2$, ... .

\begin{lm}
\label{v0R}
For any $v\in K[M]$ there exist $w_1\in K[M]$, $w_2,w_3\in K[M_1]$ such that $v=v_0w_1+x_0w_2+w_3$.
\end{lm}

\begin{prf}
The idea of the proof is analogous to The Remainder Theorem. We are going to show that $v=x_0x_2+w_3$ modulo the principal right ideal $v_0R$, where $R=K[M]$. All equalities below will be done modulo this ideal. Each of them can be multiplied on the right.

We know that $x_0^2=\beta^{-1}(\alpha x_0x_1+\alpha^2x_3^2+\alpha^2\beta x_3x_4)$. That is, $x_0^2=x_0\xi_2+\eta_2$ for some $\xi_2,\eta_2\in K[M_1]$. Let us prove by induction on $k\ge2$ that $x_0^k=x_0\xi_k+\eta_k$ for some $\xi_k,\eta_k\in K[M_1]$. We know that for $k=2$. Let this equality hold for some $k\ge2$. Then $x_0^{k+1}=(x_0\xi_k+\eta_k)x_0=x_0^2\phi(\xi_k)+x_0\phi(\eta_k)$, where $\phi$ is an endomorphism that takes each $x_i$ to $x_{i+1}$ ($i\ge0$). Therefore, $x_0^{k+1}=(x_0\xi_2+\eta_2)\phi(\xi_k)+x_0\phi(\eta_k)=x_0\xi_{k+1}+\eta_{k+1}$, where $\xi_{k+1}=\xi_2\phi(\xi_k)+\phi(\eta_k)$, $\eta_{k+1}=\eta_2\phi(\xi_k)$. 

Now we can decompose $v$ by powers of $x_0$, that is, $v=\zeta_0+x_0\zeta_1+x_0^2\zeta_2+\cdots+x_0^d\zeta_d$ for some $d$, where $\zeta_0,\ldots,\zeta_d\in K[M_1]$. Replacing here all $x_0^k$ by $x_0\xi_k+\eta_k$ for $k\ge2$, we get the desrired equality of the form $v=x_0w_2+w_3$ modulo the ideal $v_0R$, where $w_2,w_3\in K[M_1]$. Hence there exists $w_1\in R$ such that $v=v_0w_1+x_0w_2+w_3$.  

The proof is complete.
\end{prf}
\vspace{1ex}

Let $(u,v)$ be any solution of~(\ref{012}). We have $au=bv=b(v_0w_1+x_0w_2+w_3)$ according to Lemma~\ref{v0R}, together with $au_0=bv_0$ multiplied by $w_1$ on the right. This implies $a(u-u_0w_1)=b(x_0w_2+w_3)$. Recall that $b=x_1+\beta x_2$, so the right-hand side belongs to $R_1+x_0R_1$, where $R_1=K[M_1]$. Since $a=x_0+\alpha x_2$, the term $w=u-u_0w_1$ in its normal form cannot involve $x_0$. So it belongs to $R_1$. 

We have $(x_0+\alpha x_2)w=(x_1+\beta x_2)(x_0w_2+w_3)$, that is, $x_0w+\alpha x_2w=x_0(x_2+\beta x_3)w_2+(x_1+\beta x_2)w_3$. 
Comparing coefficients, we get $w=(x_2+\beta x_3)w_2$ and $\alpha x_2w=(x_1+\beta x_2)w_3$. Hence $\alpha x_2(x_2+\beta x_3)w_2=(x_1+\beta x_2)w_3$. This implies that $x_1w_3$ is left divisible by $x_2$. This is also true for all monomials involved in $x_1w_3$. From elementary properties of the monoid $M$ it follows that $w_3$ is left divisible by $x_3$, so we can put $w_3=x_3w_4$ for some $w_4\in K[M_1]$. Now $x_1w_3=x_1x_3w_4=x_2x_1w_4$, and we can cancel the equality by $x_2$ on the left. 

Now we get $\alpha(x_2+\beta x_3)w_2=(x_1+\beta x_3)w_4$. Notice that $u$, $v$ are homogeneous of some degree $d$. So $\deg x_0w_2=\deg v=d$. Therefore $\deg w_2=d-1$. Also $\deg w_4=\deg w_3-1=\deg v-1=d-1$. So we are trying to solve the equation $(x_1+\beta x_3)w_4=(x_2+\beta x_3)(\alpha w_2)$ in $K[M_1]$ in polynomials of degree $d-1$. One can apply $\phi^{-1}$ decreasing all subscripts by 1. This gives us equation $(x_0+\beta x_2)\phi^{-1}(w_4)=(x_1+\beta x_2)\phi^{-1}(\alpha w_2)$ in $K[M]$ for which we eventually know all solutions of degree $<d$ by induction. Any such solution brings us a solution $(u,v)$ of the original equation of degree $d$.
\vspace{1ex}

Let ${\rm Sol\,}(d,\alpha,\beta)$ denote the set of solutions of~(\ref{012}) with parameters $\alpha$, $\beta$, where $\deg u=\deg v=d$. This set is zero for $d=0$ and $d=1$. Let $(u',v')$ be arbitrtary solution from ${\rm Sol\,}(d-1,\beta,\beta)$. We assume by induction this set is known for us. Now we can express $w_4=\phi(u')$ and $\alpha w_2=\phi(v')$. From $\alpha x_2w=(x_1+\beta x_2)w_3=(x_1+\beta x_2)x_3w_4=x_2(x_1+\beta x_3)w_4$, we get $w=\alpha^{-1}(x_1+\beta x_3)\phi(u')$. Therefore, $u=u_0w_1+w=u_0w_1+\alpha^{-1}(x_1+\beta x_3)\phi(u')$ and $v=v_0w_1+x_0w_2+w_3=v_0w_1+\alpha^{-1}x_0\phi(v')+x_3\phi(u')$. These equalities 
\be{eqs}
\left\{
\begin{array}{ccl}
u	& = & u_0w_1+\alpha^{-1}(x_1+\beta x_3)\phi(u') \\
v	& = & v_0w_1+\alpha^{-1}x_0\phi(v')+x_3\phi(u')
\end{array}
\right.
\ee
allow us to describe all solutions of degree $d$ for parameters $\alpha$, $\beta$ provided all solutions of degree $d-1$ for parameters $\beta$, $\beta$ are known.

We see that the general case of parameters is reduced to the case $\alpha=\beta$. Also if we know the first component $u$ of the solution, then the second component $v$ can be uniquely express from it. To present information in a compact from, we formulate the following statement.

\begin{thm}
\label{descr}
Let $\beta\ne0$ be an element of a field $K$. Let us consider the equation $(x_0+\beta x_2)u=(x_1+\beta x_2)v$ in the monoid ring $K[M]$. Let
$$
u_0=x_0x_3+\beta x_0x_4-x_1x_3-\beta x_1x_4-\beta x_3^2-\beta^2 x_3x_4,
$$
$$
v_0=x_0^2-x_0x_1-\beta x_3^2-\beta^2x_3x_4
$$
be its basic solution.

Then for any its solution, one has the following presentation for its first unknown:
$$
u=u_0w_0+(x_1+\beta x_3)\phi(u_0)w_1+(x_1+\beta x_3)(x_2+\beta x_4)\phi^2(u_0)w_2+\cdots+\prod\limits_{i=1}^k(x_i+\beta x_{i+2})\phi^k(u_0)w_k,
$$
where $k\ge0$, and $w_i$ belongs to $K[M_i]$, where $M_i$ is the submonoid of $M$ generated by $x_i,x_{i+1},\ldots\,$ $(0\le i\le k)$.
\end{thm}

This result immediately follows by induction from~(\ref{eqs}). Notice that we cancelled the basic solution by a nonzero constant. Also the constant multiple in~(\ref{eqs}) does not play any role since it can be sent to $w_i$.

\end{document}